\documentclass[11pt]{amsart}

\usepackage{xcolor}
\usepackage{amsmath, amsthm, amssymb, amsfonts, mathtools, enumerate}
\usepackage{hyperref}
\usepackage{verbatim} 
\usepackage{graphicx,epsfig,color}
\usepackage{exscale}

\theoremstyle{plain}

\newtheorem{theorem}{Theorem}[section]
\newtheorem*{th:re}{Theorem \ref{th:re}}

\newtheorem{corollary}[theorem]{Corollary}
\newtheorem{conjecture}[theorem]{Conjecture}
\newtheorem{proposition}[theorem]{Proposition}

\theoremstyle{definition}

\newtheorem{definition}[theorem]{Definition}

\theoremstyle{remark}
\newtheorem{remark}[theorem]{Remark}

\newtheorem{example}[theorem]{Example}
\newtheorem{examples}[theorem]{Examples}

\DeclareSymbolFont{AMSb}{U}{msb}{m}{n}
\DeclareMathSymbol{\N}{\mathalpha}{AMSb}{"4E}
\DeclareMathSymbol{\R}{\mathalpha}{AMSb}{"52}
\DeclareMathSymbol{\Z}{\mathalpha}{AMSb}{"5A}
\DeclareMathSymbol{\D}{\mathalpha}{AMSb}{"44}
\DeclareMathSymbol{\s}{\mathalpha}{AMSb}{"53}

\renewcommand{\Im}{\mbox{Im}}

\newcommand{\sF}{{\scriptscriptstyle{F}}}

\newcommand{\sB}{{\scriptscriptstyle B}}
\newcommand{\sX}{{\scriptscriptstyle{X}}}
\newcommand{\sM}{{\scriptscriptstyle{M}}}
\newcommand{\sN}{{\scriptscriptstyle N}}

\DeclareMathOperator{\nor}{nor}

\DeclareMathOperator{\vol}{vol}

\DeclareMathOperator{\id}{id}

\DeclareMathOperator{\supp}{supp}

\DeclareMathOperator{\de}{d}

\DeclareMathOperator{\m}{m}
\DeclareMathOperator{\ric}{ric}

\newcommand{\mwp}{B\times_f F}
\newcommand{\nwp}{B\times_f^\sN F}
\newcommand{\smwp}{{\scriptscriptstyle{B\times_f F}}}

\newcommand{\dint}{\de}

\newcommand{\RCD}{\mathsf{RCD}}
\newcommand{\MCP}{\mathsf{MCP}}
\newcommand{\CD}{\mathsf{CD}}

\newcommand{\Geo}{\mbox{Geo}}
\title[curvature bounds for warped products]{Warped products and synthetic lower curvature bounds: an overview}
\author{Christian Ketterer}
\thanks{ This is an invited review for the conference “Metric measure spaces, Ricci curvature, and optimal
transport” at Villa Monastero, Varenna, operated from September 23, 2024 to October 1, 2024.
I sincerely thank Fabio Cavalletti, Matthias Erbar, Jan Maas, and Karl-Theodor Sturm for their
organization of this stimulating event and the great opportunity to compose this survey. 
This overview was
completed during the HIM Trimester Program “Metric Analysis” at the Hausdorff Research In-
stitute for Mathematics, Bonn. I sincerely thank the institute for creating such a stimulating
research atmosphere.}
\thanks{{\it 2020 Mathmatics Subject Classification.} Primary  49Q22, 53C21, 83C99. {\it Keywords}: warped products, metric spaces, curvature bounds, Lorentzian geometry}
\address{}
\email{christian.ketterer@mu.ie}
\begin{document}
\begin{abstract} This is a survey about the construction of warped products  between (semi-)Riemannian manifolds and metric (measure) spaces. The resulting  spaces  will be semi-Riemannian manifolds,  metric (measure) spaces or Lorentzian metric and metric measure  spaces. We present details of the contruction in each case and we will highlight important properties like fiber independence and the energy equation. Warped products behave nicely in relation with  curvature lower bounds. Here  we will focus on sectional and Ricci curvature lower bounds and their Lorentzian counterparts. Throughout the article we provide many examples and formulate questions and conjectures. 
\end{abstract}
\maketitle
\tableofcontents
Warped products  are a natural  generalization  of 
products of pairs of spaces. For instance, given two Riemannian  manifolds $B$ and $F$, the base and the fiber, a warped product of $B$ and $F$ is obtained by the product space $B\times F$ equipped with a Riemannian metric that arises from the product metric by scaling the Riemannian metric $g_\sF$ on each fiber $\{r\}\times F$ according to a factor given by $f(r)>0$ for a smooth function $f$ on $B$. 
%
%

Warped products play an important role in Riemannian geometry. They give a wide class of examples of spaces with lower and upper curvature bounds, including sectional, Ricci and Scalar curvature. They arise as special  spaces in rigidity and almost rigidity theorems under lower curvature bounds.  
And they are  a particular case of a Riemannian submersion and  can serve as building blocks in gluing constructions, e.g. \cite{wraithsurg}. 

{\color{black} Euclidean cones and suspensions are special cases of warped products of metric spaces that were known to respect   curvature bounds in the sense of Alexandrov \cite{bbi, bgp}.
Famously Cheeger and Colding have shown that  Riemannian
manifolds satisfying certain lower Ricci curvature bounds and almost extremal volume
or diameter conditions are Gromov–Hausdorff close to warped products of
intervals with intrinsic metric spaces \cite{almostrigidity, cheegercoldingI, coldingshape}.  Chen and Alexander-Bishop then studied general warped products of  intrinsic metric spaces and their relation to  Alexandrov curvature bounds \cite{chenwp, albi0, albi, conesplittingtheorems}.

The warped product construction also works for semi-Riemannian manifolds and in the case of a $1$-dimensional base one obtains a generalized cone that is a Lorentzian manifold. Prominent examples of warped products that satisfy the Einstein equation are the Friedmann-Lemaître-Robertson-Walker cosmological models or the Schwarzschild solution that models the gravitational field of a star \cite{oneillsemi, schwarzschild}.
Finally the construction of warped products  between metric and metric measure spaces  yields new examples of non-smooth  Lorentzian  length spaces in the sense of Kunzinger and S\"amann \cite{kusae}.

{\color{black} In this paper we will survey the construction of warped products in smooth and non-smooth context, present  examples, explain crucial properties of warped products and  highlight the connection to synthetic  sectional  and Ricci curvature lower bounds. Recently warped products also appeared in the context of lower Scalar curvature bounds (for instance \cite{gromovgafa, cz}) but this will not be addressed in this article.}

In Section 1 we recall the construction for smooth spaces, in Section 2 we explain in detail the construction for intrinsic metric spaces resulting in  nonsmooth metric or metric measure spaces, and in Section 3 we will see how a modification of the contruction in Section 3 yields a Lorentzian length space in the sense of Kunzinger and S\"amann.
\section{Semi-Riemannian manifolds}
\noindent
The following basic properties of warped products are taken from \cite{oneillsemi}.

Given  two semi-Riemannian manifolds $B$ and $F$ a  $2$-tensor on the product manifold $B\times F$ can be defined via 
$ P_\sB^* g_\sB + P_\sF^*g_\sF$
where $P_\sB$ and $P_\sF$ are the projections of $B\times F$ onto $B$ and $F$, respectively.

We use the notation $\langle \cdot, \cdot\rangle_\sB=\langle \cdot, \cdot \rangle$ for the metric $g_\sB$ on $B$, and the notation  $\langle \cdot, \cdot\rangle_\sF=(\cdot, \cdot)$ for the metric $g_\sF$ on $F$, as well as $D^\sB$ and $\nabla$ for the Levi-Civita connections on $B$ and $F$, respectively.

 The manifold $B\times F$ equipped with the semi-Riemannian metric tensor 
$$g_{\smwp}:= P_\sB^* g_\sB + (f\circ P_\sB)^2 P_\sF^* g_\sF$$ 
is called semi-Riemannian warped product $B\times_f F$ between $B, f$ and $F$. If $v, w$ are tangent vectors of $B\times F$ at $(p, x)$, then 
$$\langle v, w\rangle := g_\smwp(v,w)= \langle  DP_\sB v, DP_\sB w\rangle  + (f(p))^2  (DP_\sF v, DP_\sF w).$$ 
The Riemannian manifold $B$ is called the base and $F$ is called  fiber space or cross section. 
Clearly, the warped product $B\times_fF$ is a Riemannian manifold if $B$ and $F$ are Riemannian manifolds.

The embedded fibers $\{p\}\times F$ and the embedded leaves $B\times \{x\}$ are semi-Riemannian submanifolds of $B\times_f F$. The warped product metric $g_\smwp$ is characterized by 
\begin{enumerate}
\item $P_\sB|_{B\times \{x\}}$ is a Riemannian isometry onto $B$ for every $x\in F$. 
\item $P_\sF|_{\{p\}\times F}$ is a positive homothety map onto $F$ with scale factor $\frac{1}{f(p)}$ for every $p\in B$. 
\item For each $(p,x)\in B\times F$, the leaf $B\times\{q\}$ and the fiber $\{p\}\times F$ are orthogonal at $(p,q)$. 
\end{enumerate}

Vectors tangent to leaves are called horizontal, vectors tangent to fibers are called vertical. The horizontal  (or vertical) projection of tangent vectors $v\in T(B\times_f F)$ to $TB$ (or $TF$) is denoted with $\nor(v)$ (or $\tan(v)$).  If $V$ and $W$ are vertical vector fields $T_{(r,x)}(B\times F)$, then $II(V,W)= \nor(D_V W)$ is the shape operator (or the second fundamental form of the fiber that contains $(r,x)$)  where $D$ is the Levi-Civita connection of $B\times_f F$. 

There is a notion of lifting   a vector field on $B$ or on $F$ to a vector field on $B\times F$. The set of all such lifts is denoted with $\mathcal L(B)$ and $\mathcal L(F)$, respectively. 
If $h\in C^1(B)$, then  the vector field $\nabla (h\circ P_\sB)$ is the lift of $\nabla h$ to $B\times F$. 

Let $D$ be  the Levi-Civita connection of $B\times_f F$. If $X, Y\in \mathcal L(B)$ and $V,W\in \mathcal L(F)$ it  follows by a straightforward computation
\begin{enumerate}
\item $D_XY\in \mathcal L(B)$ is the the lift of  $D^\sB_XY$ on $B$,
\item $D_XV= D_VX = X(\log f) V$,
\item $\nor D_VW= II(V,W)= - \langle V, W\rangle \nabla \log f$,
\item $\tan D_VW\in \mathcal L(F)$ is the lift of $\nabla_V W$ on $F$. 
\end{enumerate}

As corollary we get that the leaves $B\times \{x\}$ are totally geodesic. In particular, $B$ embeds via a distance preserving map into $B\times \{x\}$.  Hence, we can write  $D= D^\sB$ in the following. On the other hand, the fibers $\{p\}\times F$ are totally umbilic. 

\begin{examples}\label{examples1}\hspace{0mm}
(1) A surface of revolution is a warped product. More precisely, if a surface $S\subset \R^3$ is obtained by revolving the image  $C=\Im \gamma$ of a smooth curve $\gamma:I\rightarrow \R^3$ in the $xz$-plane around the $x$-axis, and $f(t)$ is the distance of $\gamma(t)$ to the $x$-axis, then $S\simeq I\times_f \mathbb S^1$ where $\mathbb S^1$ is the circle of radius $1$. 
%
\smallskip\\
(2) The punctered Euclidean space $\R^3\backslash \{0\}$ is isometric to the  warped product $ (0, \infty)\times_{\id} \mathbb S^2$ where $\id(r)=r$. Spherical coordinates on $\R^3\backslash \{0\}$ are given through the map 
$$(r, t, \phi)\mapsto (r \sin t, r\sin\phi \cos t, r\cos \phi \cos t).$$
In these coordinates the standard Euclidean metric  takes the form of  a warped product metric
$$(dr)^2 + r^2 \left( (d t)^2  + \sin^2(t)( d \phi)^2\right).$$
Alternatively, since $\mathbb S^2 \simeq \{x\in \R^3: |x|_2=1\},$ we can identify  $\R^3\backslash \{0\}$ with $(0, \infty)\times \mathbb S^2$ using  the map $(r,x)\mapsto rx$. The pull-back of the Euclidean metric under this map is the warped product metric $(dr)^2 + r^2 g_{\mathbb S^2}.$
\smallskip\\
(3) The subset $\{x\in \R^3: \langle x, x \rangle < 0\}=: \mathbb M^3$ in the punctured Minkowski space $\R^3_1\backslash \{0\}= (\R^3\backslash \{0\}, \langle \cdot, \cdot \rangle_1)$ is a warped product where $\langle v, w\rangle_1= -v_1 w_1 + v_2 w_2 + v_3 w_3$.  We can identity the $2$-dimensional hyperbolic plane with  $\{ x\in \R^3_1: \langle x, x \rangle_1=-1\}= \mathbb H^2\subset \R^3_1$ furnished with the restriction of $\langle \cdot, \cdot \rangle_1$.  Considering $(r,x) \in (0,\infty)\times \mathbb H^2\mapsto rx$ the pull-back of $\langle \cdot , \cdot \rangle_1$ under this map takes the form of warped product
$(dr)^2 + r^2 g_{\mathbb H^2}.$
\smallskip\\
(4) Let $\mathbb S^2_{\geq 0}= \mathbb S^2 \cap \{(x, y, z)\in \R^3: x\geq 0\}$ where we see $\mathbb S^2$ as the unit sphere w.r.t. the Euclidean norm on $\R^3$. Let $\partial \mathbb S^2_{\geq 0} = \mathbb S^2\cap \{(x,y,z)\in \R^3: x=0\}$ be the boundary  of $\mathbb S^2_{\geq 0}$, and let $d_{\partial \mathbb S^2_{\geq 0}}$ be the distance function in $\mathbb S^2_{\geq 0}$ to the boundary. 
If we parametrize $\mathbb S^2_{\geq 0}\backslash \{(1,0,0)\}$ via $\Phi: (t,\phi)\in [0,\frac{\pi}{2})\times \R\mapsto ( \cos t, \sin \phi \sin t, \cos \phi \sin t)$ then $d_{\partial \mathbb S^2_{\geq 0}}(x)=t=t(x)$. If we set $f(x) = \sin d_{\partial \mathbb S^2_{\geq 0}}(x)$,  $f$ is a smooth function that is positive on $\mathbb S^+_{>0}$. Then the warped product $\mathbb S^2_{>0}\times_f \mathbb (0,\pi)$ is isometric to  $\mathbb S^3_{>0} \backslash \{(1,0,0, 0)\}\subset \R^4$ where $\mathbb S^3_{> 0} = \{(x, y, z, w): x>0\}$.
\end{examples}
\subsection{Warped product geodesics}
A curve $\gamma=(\alpha, \beta)$ in $M=B\times_f F$ is a geodesic, i.e. $D_{\gamma'}\gamma'=0$, if and only if
\begin{enumerate}
\item $\alpha'' = (\beta' , \beta') (f\circ \alpha) \nabla f\circ \alpha$, 
\item $\beta''= \frac{-2}{f\circ \alpha}  \frac{d(f\circ \alpha)}{ds} \beta'$. 
\end{enumerate}
In particular, $$\frac{d}{ds} (f\circ \alpha)^4 (\beta', \beta')= 4 (f\circ \alpha)^3 \frac{d( f\circ \alpha)  }{ds}(\beta', \beta') + (f\circ \alpha)^4 2 (\beta'', \beta)=0.$$ Hence $(f\circ \alpha)^4 (\beta', \beta')= (f\circ \alpha)^2 \langle \beta', \beta'\rangle$ is a constant $C$, i.e. the speed of $\beta$ is proportional to $\frac{1}{2f^2\circ \alpha}$. 
\begin{remark}\label{rem:geod} The second property  above implies that $\beta$ is a pregeodesic in $F$ (i.e. a geodesic up to reparametrization).
It also follows that 
$$\alpha'' =\frac{C}{(f\circ \alpha)^3} \nabla f\circ \alpha=- \nabla \left( \frac{C}{2f^2}\right)\circ \alpha,$$
i.e. $\alpha$ is a gradient flow curve of $\frac{C}{2f^2}$, and in partiular it is independent of $\beta$ up to reparametrization.  After an affine reparametrization of $\gamma$ we can set $C=+1, 0, -1$. Moreover, since $\gamma$ is still  a geodesic, it follows 
$$\frac{1}{2} \langle \alpha', \alpha'\rangle + \frac{1}{2 (f\circ \alpha)^2} =E.$$The sign of $C=+1, 0, -1$ is determined by the sign of $(\beta', \beta')$,  the causal character of $\beta$. In particular, if $F$ is a Riemannian manifold, then either $C>0$ or $C=0$. In the latter case $\beta$ is constant and $\alpha$ is a geodesic in $B$. 
\end{remark}
\begin{remark}\label{rem:clairaut}
If we consider a surfaces of revolution and $\gamma$ such that $\beta'>0$, then the identity $$(f\circ \alpha)(t) \beta'(t)=const$$ is known as  {\it Clairaut's relation}. Indeed, with the notation we used before $f\circ \alpha(t) \beta'(t)$ is the inner product of $\gamma'(t)$ with ${\partial_\phi}(\gamma(t))$. Hence $$const = \sqrt {\langle \partial_\phi(\gamma(t)), \partial_\phi(\gamma(t))\rangle}\cos \theta(t)= f(\gamma(t))\cos\theta(t)$$ where $\theta(t)$ is the angle between $\gamma'(t)$ and the meridian through $\gamma(t)$. 
\end{remark}
%

For the rest of the exposition we will assume the following.  Let $B$ and $F$ be both Riemannian manifolds with Riemannian metrics $g_B$ and $g_F$. We  can consider $B$ as a semi-Riemannian manifold  that  is either equipped with the Riemannian metric $g_B$ or with the semi-Riemannian $-g_\sB= g_{-\sB}$.  We denote these spaces with $B=+B$ and $-B$. The warped product metric  $\pm B\times_f F$   is $B\times F$ equipped with the metric
$${\textstyle g_{\pm \smwp} = (P_\sB)^* g_{^\pm \sB} + (f\circ P_\sB)^2 (P_\sF)^* g_\sF= \pm (P_\sB)^* g_{\sB} + (f\circ P_\sB)^2 (P_\sF)^* g_\sF.}$$
In the following we write $\pm B \times_f F$ to describe both cases where $B$ is equipped with $g_\sB$ or with $g_{-\sB}$. 

\subsection{Curvature of warped products} 
\subsubsection{Sectional curvature}We recall the definition of the curvature tensor. If $(M, g_M)$ is a semi-Riemannian manifold, then the curvature tensor is defined as 
$$R(X,Y)Z= \nabla_X\nabla_Y Z - \nabla_Y \nabla_X Z- \nabla_{[X,Y]}Z$$
for vector fields $X,Y,Z$ and the Levi-Civita connection $\nabla$ of $M$.

The sectional curvature of a $2$-plane $P$ in $TM$ is  
$$K(P)= \frac{g_{\sM}(R(v,w)v,w)}{g_\sM(v,v) g_\sM(w,w)}$$
where $v,w$ are orthogonal vectors that span $P$.  

We say $M$ has sectional curvature is bounded from below by $K$ if 
\begin{align}\label{secbd} g_{\sM}(R(v,w)v,w)\geq K g_{\sM}(v,v) g_{\sM}(w,w) \ \forall v\perp w.\end{align}
\begin{remark}
If $M$ is a Riemannian manifold, then $M$ has sectional curvature bounded from below if and only if $K(P)\geq K$ for all $2$-planes $P$ in $TM$. 

However, if $M$ is a semi-Riemannian manifold, then $M$ has sectional curvature bounded from below by $K$ if and only if 
\begin{align*}
K(P)\begin{cases} \geq K & \mbox{ if } P \mbox{ is spacelike}; \\
\leq K & \mbox{ if } P \mbox{ is non-degenerated and not space-like}. 
\end{cases}
\end{align*}
A $2$-plane is called spacelike like if it is spanned by  vectors $v$ and $w$ with positive length w.r.t. $g_{\sM}$. It is nondegenerated if it spanned by two $v$ and $w$ that both have non-zero length w.r.t. $g_{\sM}$. 

We say that a Lorentzian manifold $M$ of signature $(-, +, \dots, +)$  has timelike sectional curvature bounded from below by $K$ if \eqref{secbd} holds for all $v\perp w$ such that $v,w$ span a nondegenerated, non-space-like plane. 
\end{remark}
We consider the warped product $\pm B\times_f F$ where $B$ and $F$ are Riemannian manifolds and $f>0$ is smooth function on $B$. The sectional curvature $K_\smwp=K$ of $\pm \mwp$ was displayed by Bishop and O'Neill in \cite{bio}. 
\begin{proposition}[Bishop-O'Neill] Let
$X,Y\in T_rB$ with $\langle X, Y\rangle_B=0$,  and $v, w\in T_xF$ with $\langle v,w\rangle_F=0$. 

We  assume
$\langle (X,v), (Y,w)\rangle_{\pm B\times_f F}=0$ and consider the $2$-plane
$P= \mbox{Span} ((X,v), (Y,w))\subset T_{(r,x)} (\pm B\times_f F)$. Then the sectional curvature $K(\cdot)$ of the warped product $\pm B\times_f F$ satisfies
\begin{align*}K(P)=& \pm K_B(X,Y) |X|_B^2|Y|_B^2 \\
&  - f(r) \Big[ |w|_F^2 \nabla^2 f(X,X) + |v|_F^2 \nabla^2 f(Y,Y)\Big] \\
&  + f^2(r)\Big[K_F(v,w)\pm |\nabla f|^2_B(r)\Big] |v|^2_F |w|_F^2.
\end{align*}
\end{proposition}
In \cite{albi,albilo} Bishop and Alexander deduce the following corollary that is a straightforward consequence of the previous proposition.
\begin{corollary}[Bishop-Alexander]
The warped product 
$\pm B\times_f F$ has sectional curvature bounded from below by $K$ if and only if 
\begin{enumerate}
\item $\nabla^2 f \pm K f g_B\leq 0$
\item $K_F(\cdot) \geq Kf(r)^2 \pm  |\nabla f|^2_B(r)$,
\item $\dim_B=1$ or $\mbox{sec}_B\geq K \ (\leq  K)$.
\end{enumerate}
\end{corollary}
\subsubsection{Ricci curvature}
If $M$ is a semi-Riemannian manifold, then
the Ricci curvature $\ric_\sM(v,v)$ in direction of a tangent vector $v$ is the mean value of the sectional curvature of planes that contain $v$. $\ric$ is $2$-tensor on $M$.  

If $M$ is a Riemannian manifold, we say that $M$ has Ricci curvature bounded from below by $K$, denoted $\ric_M\geq K$, if 
$$\ric_{\sM}(v,v) \geq K|v|_{\sM}^2 \mbox{ for all }v\in TM.$$

If $M$ is a Lorentzian manifold, we say that $M$ has timelike Ricci curvature bounded from below by $K$ if 
$$\ric_{\sM}(v,v)\geq K |v|^2_{\sM}\mbox{ for all } v\in TM \mbox{ s.t. } g_{\sM}(v,v)<0$$where $|v|_M= \sqrt{-g_M(v,v)}$.
\begin{remark}
Timelike sectional curvature bounded from below by $K$ implies timelike Ricci curvature bounded from below $-K$. 
\end{remark}

\begin{proposition}[\cite{oneillsemi}]  Let $\pm B\times_fF$ be a warped product between $\pm B$, $f$ and $F$. Then
\begin{align*}
&\ric_{B\times_f F}(X+V, X+V) = \ric_B(X,X) + n { \frac{\nabla^2 f(X,X)}{f}}\\
&\ \ \ \ \  \ \ \ \ \ + \ric_F(V,V) - \left(\frac{\Delta^{\pm B} f }{f} + (n-1)\frac{|\nabla f|^2_{\pm B}}{f^2} \right) |V|_{\pm B\times_f F}.
\end{align*}
\end{proposition}
\begin{remark}
We note that $\nabla^\sB= \nabla^{-\sB}$ and $\nabla_\sB^2 f= \nabla_{-\sB}^2 f$, i.e. the Levi-Civita connection and the Hessian operator don't change when replacing $B$ with $-B$. But it holds
\begin{center}$|\nabla f|_{\boldsymbol{\pm} B} = \boldsymbol{\pm} |\nabla f|_\sB, \ \ \Delta^{\boldsymbol{\pm} B} f=\boldsymbol{ \pm} \Delta^B f.$\end{center}
\end{remark}
\begin{corollary}
$\ric_{\pm B\times_f F} \geq \pm nK$ if 
\begin{enumerate} 
\item $ \nabla^2 f + K f g_{\pm \sB} \leq 0$,
\item $\ric_F\geq (n-1)(Kf^2  + |\nabla f|^2_{\pm B})$, 
\item $\ric_\sB\geq (d-1)K.$
\end{enumerate}
\end{corollary}
\begin{example}
One can ask if, as for sectional curvature, the conditions (1), (2) and (3) are also necessary for the the lower Ricci curvature bound of $\pm B\times_f F$. However, this is not true for Riemannian warped products as the following example shows. 
Let $B=\R$ and let $F$ be an $n$-dimensional  Riemannian manifold of constant  curvature $-1$. In particular, $F$ is Einstein with  Ricci curvature  equal to $-(n-1)$.  The  Riemannian product $\R\times F$ is  a warped product w.r.t. $f(r)\equiv 1$, and it satisfies sharp bound $\ric_{\R\times F}\geq -(n-1)=  n K$ with $K:=-\frac{n-1}{n}$.  The function $f\equiv 1$ satisfies $f'' +K f\leq 0$, and $B=\R$ satisfies $\ric_\sB=0\geq (d-1)K$. But $F$ doesn't have Ricci curvature bigger than $(n-1)\left(Kf^2 + (f')^2\right)=(n-1)K$, since $(n-1)K=-(n-1)\frac{n-1}{n} >-(n-1)$ and since $F$ is Einstein with $\ric_\sF= - (n-1) g_\sF$.
\end{example}
\begin{remark}
If one assumes $\ric_\sB\geq (d-1) K$ and $$\nabla^2 f+ K fg_{\pm \sB}=0,$$ the condition $\ric_{\pm \smwp}\geq nK$ implies 
$\ric_\sF\geq (n-1)(Kf^2 + |\nabla f|^2_{\pm B})$. For instance, this holds in the case of an Euclidean cone. 
\end{remark}
\subsubsection{Bakry-Emery Ricci curvature}

Let $M$ be a  Riemannian manifold and assume $M$ is equipped with  a measure  $\Phi \de \vol_M$ where $\Phi\in C^{\infty}(M)$ and $\Phi> 0$.  We call   $(M, \Phi)$ a weighted Riemannian manifold. 

Let $N\geq 1$.
If  $N>n=\dim_M$,  for $p\in  M$ and $v\in T_pM$ the Bakry-Emery $N$-Ricci tensor is defined through
\begin{align*}
\ric_g^{\Phi, N}|_p(v,v)
= \ric_\sM|_p(v,v) - (N-n) \frac{ \nabla^2 \Phi^{\frac{1}{N-n}}|_p(v,v)}{\Phi^{\frac{1}{N-n}}(p)}.
\end{align*}
If $N=n$, we set 
\begin{align*}
\ric_\sM^{\Phi, n}|_p(v,v) = \begin{cases} \ric_g|_p(v,v) - \nabla^2 \log \Phi |_p(v,v)
&\mbox{
if $g( \nabla \log \Phi _p, v)=0$}, \\
-\infty&\mbox{ otherwise.}
\end{cases}\end{align*} 
Finally, let $\ric_g^{\Phi, N}\equiv-\infty$ if $1\leq N<n$. 
\subsubsection{$N$-warped products}\label{subsubsec:Nwp}
Let $B$ and $f$ be as before, and let $(F, \Phi)$ be a weighted Riemannian manifold. We can turn the warped product $\mwp$ into a weighted Riemannian manifold by introducing the measure 
$$\m^\sN= f^\sN \de\vol_\sB \Phi \de\vol_\sF= f^{N-n} \Phi \de \vol_\smwp= \Theta \de \vol_{\smwp}$$
on $B\times_f F$ where $\dim_\sF=n$ and $\Theta:= f^{N-n} \Phi $.

We call  the weighted Riemannian manifold $(\mwp, \m^\sN)=: \nwp$ the $N$-warped product of $B$, $f$ and $(F, \Phi)$.

\begin{proposition}[{\cite{ketterer1}}]\label{prop:alt} Let $\nwp$ be  an $N$-warped product. Then
\begin{align}\label{formula:bakry}
&\ric_{B\times_f F}^{\Theta, N+d}(X+V, X+V)  =%
\nonumber\\
&\ric_\sB^{f^{N-n}, N}(X,X)+ \ric_F^{\Phi, \sN}(V,V) - \left(\frac{\Delta^\sB f }{f} + (N-1)\frac{|\nabla f|^2_{ B}}{f^2} \right) |V|_{B\times_f F}
\end{align}
where we note that $\ric_\sB^{f^{N-n}, N}(X,X)=\ric_B(X,X) + N { \frac{\nabla^2 f(X,X)}{f}}.$
\\
Consequently, $\ric_{\smwp}^{\Theta, N+d}\geq (N+d-1)K$ if 
\begin{enumerate}
\item $\Delta^\sB f+d K f\leq 0 .$
\item $\ric_\sF^{\Phi, N}\geq (N-1)(Kf^2  + |\nabla f|^2_{\pm\sB})$
\item $\ric_\sB^{f^{N-n}, N}\geq (N + d-1)K$.
\end{enumerate}
\end{proposition}
\begin{remark}
The conditions (1) and (3) together are implied if \begin{center}$\ric_{\sB}\geq (d-1) K$ and $\nabla^2 f+ Kfg_\sB\leq 0$.
\end{center}
\end{remark}
Exactly as in the positive signature case one can consider weighted semi-Riemannian, or more specifically, Lorenztian manifolds and their Bakry-Emery $N$-Ricci curvature for $N\geq 1$ (or even with negative $N$). 
\begin{conjecture}
A formula  like  \eqref{formula:bakry} for the Bakry-Emery $N$-Ricci tensor holds in the context of weighted semi-Riemannian manifolds.
\end{conjecture}
\section{Metric  and metric measure spaces}
In this section we outline how to define  warped products as metric spaces between general, intrinsic metric spaces as in \cite{bbi, albi0, albi}.
\subsection{Length spaces}
To motivate the  construction  we recall first  the construction of a length structure on a metric space.

Let $X$ be a metric space and let $\gamma:[a,b]\rightarrow X$ be a continuous map. We call $\gamma$ a path in $X$. The length of $\gamma$ is defined via
$$
\sup \sum_{i=1}^N \de_\sX(\gamma(t_{i-1}), \gamma(t_i))=: L^{\sX}(\gamma)$$
where the supremum is w.r.t. to  partitions $a=t_0\leq t_1 \leq   \dots \leq t_N=b$ of $[a,b]$. 

The length functional $L^\sX$ is a lower semi-continuous length structure in the sense of \cite{bbi}.  In particular, $L^\sX(\gamma)$ is invariant under suitable reparametrizations of $\gamma$. If $L^\sX(\gamma)<\infty$, $\gamma$ is called rectifiable.

If $\gamma:[a,b]\rightarrow X$ is a Lipschitz curve, then the metric speed 
$$\lim_{h\rightarrow 0} \frac{\de_\sX(\gamma(t+h), \gamma(t))}{h}=: |\gamma'(t)|$$
exists for $\mathcal L^1$-a.e. $t\in [a,b]$, and $L^\sX(\gamma)= \int_a^b|\gamma'|(t)\de t.$ In particular, every rectifiable curve admits a reparametrization that is Lipschitz \cite[Theorem 2.7.6]{bbi}.

A metric space $X$ is called  intrinsic (or a length space) if  for every pair of points $x,y\in X$ it holds $\de_\sX(x,y)= \inf L^{\sX}(\gamma)$ where the infimum is w.r.t. all rectifiable curves
whose endpoints are $x$ and $y$.  
Assuming $x,y\in X$ admit a rectifiable curve $
\gamma$ connecting them such that $L^\sX(\gamma)=\de_\sX(x,y)$, $\gamma$ is called a minimizer of $L^\sX(\gamma)$ or minimal geodesic (in the context of metric spaces we  just say geodesic). If every pair $x,y\in X$ admits a minimal geodesic connecting them,  we call $X$ strictly intrinsic or a geodesic  space. 

By the metric Hopf-Rinow theorem (Hopf–Rinow–Cohn-Vossen Theorem) \cite[Theorem 2.5.28]{bbi} a locally compact, complete and intrinsic metric space is strictly intrinsic.

Let 
$$\mbox{Geo}(X)= \{\gamma: [0,1]\rightarrow X: \gamma \mbox{ is a constant speed geodesic}\}$$
equipped with the uniform distance and let $e_t: \mbox{Geo}(X)\rightarrow X$, $e_t(\gamma)=\gamma(t)$ be the evaluation map that is continuous. 

A subset $\Upsilon\subset \mbox{Geo}(X)$ is nonbranching if for every pair $\gamma_0, \gamma_1\in \Upsilon$ it follows that $\gamma_1(\epsilon)=\gamma_2(\epsilon)$ for some $\epsilon>0$ implies that $\gamma_0=\gamma_1$. If $\mbox{Geo}(X)$ is nonbranching, we say the metric space $X$ is nonbranching. 

\subsection{Warped products between metric spaces}
Let $B$ and $F$ be intrinsic metric spaces, and let $f: B\rightarrow [0, \infty)$ be continuous. 
Let $\gamma= (\alpha, \beta):[a,b]\rightarrow B\times F$ be a continuous curve in $B\times F$.

We define the length 
$L(\gamma)$ of $\gamma$ as
$$
\sup\sum_{i=1}^N \left( \de_\sB(\alpha(t_i), \alpha(t_{i-1}))^2 + (f\circ \alpha)^2(t_i) \de_\sF(\beta(t_i) \beta(t_{i-1})^2\right)^{\frac{1}{2}}=:L(\gamma)$$
where the supremum is w.r.t.  partitions $a=t_0\leq  \dots\leq t_N=b$ of $[a,b]$.  

If $\alpha$ and $\beta$ are rectifiable, then $\gamma$ is rectifiable, i.e. $L(\gamma)<\infty$. In the following we call such a map $\gamma$ with $\alpha$ and $\beta$  rectifiable an {\it admissible path. }

The definition of $L$ mimics the definition of the induced length $L^\sX$ on a metric space $X$, and it is not difficult to check that $L$ satisfies the following properties:
\begin{enumerate}
\item If $c\in (a,b)$, then $L(\gamma|_{[a,c]})+ L(\gamma|_{[c,b]})= L(\gamma|_{[a,b]})$ ({\it Additivity});
\item If $\gamma$ is admissible, then $c,d\in [a,b]\mapsto L(\gamma|_{[c,d]})$ is a continuous function in $c$ and $d$;
\item $L$ is lower semi-continuous on the space of admissible paths in $B\times F$ w.r.t. to pointwise convergence of admissible paths. 
\end{enumerate}  The length of an admissible curve is invariant w.r.t. suitable reparametrizations. 

\begin{remark}$L$ has all the properties of a length structure in the sense of \cite{bbi}  except for being consistent with the topology of $B\times F$ in general. The latter is because $f$ is allowed to vanish. But if $(p,x)$ is a point in $B\times F$ where $f(p)>0$, one has the following local consistency around $(p,x)$:
\begin{enumerate}
\item[(4)] There exists a neighborhood $U_{(p,x)}$ of $(p,x)$ such that the length of paths connecting $(p,x)$ with points in the complement of $U_{(p,x)}$ is separated from $0$. 
\end{enumerate}
\end{remark}

If $\alpha$ and $\beta$ are  Lipschitz continuous in $B$ and $F$, respectively, the metric speeds $|\alpha'|$ and $|\beta'|$ are defined for a.e. $t\in [a,b]$, and one can show that 
$$L(\gamma) = \int_a^b \sqrt{ |\alpha'|^2+ (f\circ \alpha)^2 |\beta'|^2} dt.$$

The warped product metric $\de_{\smwp}$ on $B\times F$ is defined as the intrinsic metric associated to the length structure $L$, i.e. for two points $(p,x)$ and $(q,y)$ we define 
$$\de_{\smwp}((p,x), (q,y)):= \inf L(\gamma)$$
where the infimum is w.r.t. all rectifiable curves $\gamma$ that connect the points $(p,x)$ and $(q,y)$.  The infimum is finite provided there are rectifiable curves between $p$ and $q$ in $B$, and between $x$ and $y$ in $F$.  $\de_{\smwp}$ is symmetric and satisfies the $\triangle$-inequality. 

\begin{definition} The warped product metric space $B\times_f F$ between $B$, $F$ and $f$ is given by 
$$( B\times F/\sim, \de_{\smwp}) \  \mbox{where } (p,x)\sim (q,y)  \Longleftrightarrow  \de_{\smwp}((p,x), (q,y))=0.$$
By definition $B\times_f F$ is an intrinsic metric space. 
\end{definition}

\begin{remark}
%
Let $[(p,x)]\in B\times F/\sim$   be a point where  $f(p)=0$. We want to understand the relationship between the topology of $B\times F/\sim$ close to $[(p,x)]$ and the length structure $L$. 
%

If $[(q,y)]\neq [(p,x)]$ such that $p\neq q$, then  an dmissible path  $\gamma=(\alpha, \beta)$ always satisfies
$$L(\gamma)\geq L^\sB(\alpha)\geq \inf_\alpha L^\sB(\alpha)>0$$
where the last infimum is w.r.t.  rectifiable curves $\alpha$ that connect $p$ and $q$. 

If $[(q,y)]\neq [(p,x)]$ such that $x\neq y$ but $p=q$, then we have for every admissible path $\gamma$ that 
$L(\gamma)\geq f(p) L^\sF(\beta)>0$
provided $f(p)=f(q)\neq 0$. Otherwise, the infimimum of $L(\gamma)$ w.r.t. all admissible paths connecting $[(p,x)]$ and $[(q,y)]$ is $0$. The latter implies $[(p,x)]=[(q,y)]$ which is a contradiction. Hence, if $[(q,y)]\neq [(p,x)]$ such that $x\neq y$ and $p=q$, it follows $f(p)>0$. 
\smallskip

These considerations imply that   $L$ is consistent with the topology of $B\times F/\sim$, and hence $L$ is a lower semi-continuous length structure on the class of admissible paths 
\cite{bbi}. We note that every admissible path is also continuous in $B\times F/\sim$.

As  a consequence we obtain that  the induced length of $\de_{\smwp}$ coincides with the length structure $L$ according to \cite[Theorem 2.4.3]{bbi}.
\end{remark}
\begin{example} The examples that were presented before in \ref{examples1} natually fit into the metric framework of warped products. 
For instance, we can consider $[0, \infty)\times_{\id} \mathbb S^2$. The smooth warped product $(0, \infty)\times_{\id} \mathbb S^2$ embeds as a metric space into $[0, \infty)\times_{\id} \mathbb S^2.  $
\end{example}
The following theorem is an important result  for the structure of warped products between metric spaces.  It generalizes the properties of geodesics in semi-Riemannian warped products  in Remark \ref{rem:geod}  and in Remark \ref{rem:clairaut} to the context of metric warped products. 
\begin{theorem}[Alexander-Bishop, {\cite{albi0}}]\label{th:albi0}
Let $\gamma=(\alpha, \beta)$ be a minimizer w.r.t $L$ in $B\times_f F$ parametrized proportional to arclength. Assume $f>0$. Then 
\begin{enumerate}
\item[(a)] $\beta$ is a minimizer in $F$; 
\item[(b)] (Fiber independence) $\alpha$ is independent of $F$, except for the total height, i.e. the length $L^\sF(\beta)$ of $\beta$. More precisely, if $\hat F$ is another strictly intrinsic metric space and $\hat \beta$ is a minimizing geodesic in $\bar F$ with the same length and speed as $\beta$, then $(\alpha, \hat \beta)$ is a minimizer in $B\times_f\hat F$. 
\item[(c)] (Energy equation, version 1) $\beta$ has speed $\frac{c_\gamma}{f^2\circ \alpha}$ for a constant $c_\gamma$;
\item[(d)] (Energy equation, version 2) 
$\alpha$ satisfies $\frac{1}{2} |\alpha'|^2 + \frac{1}{2f^2\circ \alpha}= E$ a.e. where $E$ is the proprotionality constant of the parametrization of $\gamma$. 
\end{enumerate}
\end{theorem}The first property (a) was  shown in \cite{chenwp}.
\begin{remark} If we assume that $B$ and $F$ are locally compact, complete, strictly intrinsic metric spaces, the existence of minimizing curves $\gamma=(\alpha, \beta)$ for $L$ is guaranteed by the Arzela-Ascoli theorem. In particular, one has the following corollary
\end{remark}
\begin{corollary} If $B$ and $F$ are locally compact, complete, instrinsic metric spaces, then the warped product $\mwp$ is  a locally compact, complete and  intrinsic metric space. 
\end{corollary}
\begin{corollary} 
	Assume $B$ and $f$ are smooth. Then $\mwp$ is nonbrachning if and only if $F$ is nonbranching. 
	\end{corollary}
\begin{remark}
The Theorem \label{th:ab1} implies that conservative mechanics makes sense on any intrinsic metric space.  If $V$ is a function on $B$ that is bounded from below, let $\tilde V$ be  a positive function that differs rom $V$ by a constant. Defining $f={1}/{\sqrt{2 \tilde V}}$ the projections to $B$ of length minimiziers in $B\times_f \R$, are trajectories of the potential function $V$, up to reparametrization. 
\end{remark}
\begin{examples}[Metric cones \cite{bgp}] (1)
Given a metric space $F$, the {\it Euclidean cone} over $F$ with vertex $O$ is the quotient space $C(F)=[0, \infty)\times F/\sim$ where $(s, x)\sim (t,y) \in O$ if and only if $s=t=0$. The metric on the cone is defined from the {\it Euclidean cosine formula }
$$\de_{C(F)}((s,x), (t,y))^2= s^2 + t^2 - 2 st \cos\left( \min\{ \de_\sF(x,y), \pi\}\right).$$
{\it Claim.} Let $F$ be strictly intrinsic.  Then the Euclidean cone $C(F)$ over $F$ is isometric to the warped product $[0, \infty)\times_{\id} F$. 
\smallskip\\
The claim follows from Theorem \ref{th:albi0}. Locally in a neighborhood  $U\subset [0, \infty)\times_{\id} F$ of a point $(r,z)=p$ with $r>0$ we have the following property. If $\gamma=(\alpha, \beta):[a,b]\rightarrow U$ is a minimizing geodesic with endpoints $(s,x)$ and $(t,y)$, then it follows from (b) in Theorem \ref{th:albi0} that $L(\gamma)= L(\tilde \gamma)$ where $\tilde \gamma$ is the minimal geodesic in $\R^2\backslash \{0\}= (0, \infty)\times_{\id} (0, \pi)$ between $(s, 0)$ and $(t, L^\sF(\beta))$.  By (a) in Theorem \ref{th:albi0} $\beta$ is a minimizer in $F$ and since $F$ is strictly instrinsic, we have $L^\sF(\beta)= \de_{\sF}(\beta(a), \beta(b))$. 
Hence, by the cosine rule in $\R^2\backslash \{0\}$ it follows
$$\de_{\smwp}((s,x),(t,y))^2 = L(\gamma)^2= L(\tilde \gamma)^2= s^2 + t^2 - 2st \cos L^\sF(\beta).$$
Hence, in a sufficiently small neighborhood of $p$ the distance $\de_{\smwp}$ coincides with $\de_{C(F)}$. 
Since intrinsic metrics are fully determined by their local properties \cite[Corollary 3.1.2] {bbi} it follows that the intrinsic metrics $\de_{\smwp}$ and $\de_{C(F)}$ on $(0, \infty)\times F$ coincide. Hence, the closure of $\de_{C(F)}$ on $[0,\infty)\times F/\sim$ is isometric to $[0, \infty)\times_{\id} F$. 
\end{examples}
\noindent
(2) The {\it spherical suspension} over a metric space $F$ is  the quotient space $S(X)=[0, \pi]\times F/\sim$ where $(s, x)\sim (t, y)$ if and only if  $s=t=0$ or $s=t=\pi$. The metric on the cone is defined from the spherical cosine rule and coincides with the warped product between $[0,\pi]$, $F$ and $\sin:[0,\pi]\rightarrow [0, \infty)$ provided $F$ is strictly intrinsic. 
\smallskip\\
(3) The {\it elliptic cone} over a metric space $F$ is the quotient space $[0, \infty)\times F/\sim$ where $(s, x)\sim (t,y) \in O$ if and only if $s=t=0$. The metric on the elliptic cone is defined from the {\it Hyperbolic cosine rule }. It coincides with the warped product between $[0, \infty)$, $F$ and $\sinh: [0, \infty) \rightarrow [0, \infty)$ provided $F$ is strictly intrinsic. 
\smallskip\\
(3) The {\it parabolic cone} is the warped product $\R\times_{e^x} F$ between $\R$, $e^x$, $x\in \R$, and a metric space $F$.  It is clear that the smooth warped product $\R\times_{e^x} \R$ is isometric to the hyperpolic plane $\mathbb H^2$. Hence, by Theorem \ref{th:albi0} the distance between two points $(s,x)$ and $(t,y)$ in $\R\times_{e^x} F$ is the same as the distance between the points $(s, 0)$ and $(t, \de_\sF(x,y))$ in $\mathbb H^2\simeq \R\times_{e^x}\R$ provided $F$ is strictly instrinsic. 
\smallskip\\
(4) The {\it hyperbolic cone} is the warped product $\R\times_{\cosh} F$ between $\R, \cosh$ and a metric space $F$.
\subsection{Alexandrov spaces}
Let us recall the definition of geodesic spaces with curvature bounded from below (short CBB) by  $K$, also  denoted by $CBB(K)$. More details about Alexandrov spaces and the following definitions can be found in \cite{bbi, petsem, plaut}. It is said that a locally compact, complete and intrinsic metric space $X$ is $CBB(K)$ if 
any triangle $\triangle$ of perimeter less than $2\pi/\sqrt K$ is thinner than the unique comparison triangle in in the simply connected, $2$-dimensional model space $\mathbb M_K^2$ of constant curvature $K$. We shall only consider such spaces of finite dimension, and we call this class of spaces {\it Alexandrov spaces}.

\subsubsection{$fK$-concave functions}
A continuous function $f: X\rightarrow \R$ on an intrinsic metric space $X$ is said to be $fK$-concave if the restriction of $f$ to every minimal geodesic satisfies the differential inequality $f'' + Kf\leq 0$ in the distributional sense. In particular, the restriction of an $fK$-concave function $f$ to a geodesic has all the regularity properties of a convex function. Moreover $fK$-concavity implies that $f$ is local $\lambda$-concave, hence semi-concave on $X$. 
Petrunin and Perelman \cite{pepe} showed that a semi-concave function on an Alexandrov space $X$ is Lipschitz continuous, but for more general metric spaces this might not be true anymore. 

\subsubsection{Tangent space}
On every $n$-dimensional Alexandrov space we have that the space of directions $\Sigma_p$ at $p$ is a compact $(n-1)$-dimensional  Alexandrov with CBB by $1$. The Euclidean cone over $\Sigma_p$ coincides with the tangent cone $T_pX$ at $p$ that we obtain by blowing up the space around $p$ and extracting a Gromov-Hausdorff converging sequence. Every $v\in \Sigma_p$ is represented by a geodesic $\gamma:[0,\epsilon]\rightarrow X$ with $\gamma(0)=p$. 
 
 \subsubsection{Alexandrov boundary}
The boundary of an Alexandrov space $X$ is defined by induction over the dimension: A point $p$ in $X$ is a boundary point if the space of directions $\Sigma_p$ admits boundary points. For the start of the induction we note that a $1$-dimensional Alexandrov   is either isometric to $\mathbb S^1$, $\R$, $[a,b]$ or to $[0,\infty)$.

\subsubsection{Alexandrov differential}
Given a semi-concave function $f$ on $X$ we regard $Df_p$ as a function on the spaces of directions. $Df_p(v)$ is the derivative of $f\circ \gamma$ at $0$ for a geodesic $\gamma$ representing $v$. Then $|Df_p|$ is definded as the maximum of $Df_p$ on $\Sigma_p$. We note that $|Df_p|$ coincides with the absolute value of the gradient if $f$ is a smooth function and $X$ is a Riemannian manifold. 
\smallskip

The next theorem generalizes the statement on sectional curvature bounds for warped products in the case of {\color{black} signature $(+,\dots, +)$} to the context of Alexandrov spaces with curvature bounded from below. 

\begin{theorem}[Alexander-Bishop]\label{th:albi1} Let $B$ and $F$ be complete, locally compact intrinsic metric spaces, and let $f$ be a Lipschitz function on $B$. \smallskip\\
Then the  warped product 
$B\times_ f F$ has CBB by $K$ if and only if 
\begin{enumerate}
	\item \begin{itemize}\item[(a)] $B$ has CBB by $K$,\item[(b)] $f$ is $Kf$-concave,
\item[(c)] If $B^\dagger$ is the result of gluing two copies of $B$ along the closure of the set of boundary points where $f$ is nonvanishing, and $f^\dagger: B^\dagger \rightarrow [0,\infty)$ is the tautological extension of $f$, then $B^\dagger$ has CBB by $K$ and $f^\dagger$ is $fK$-concave. 
\end{itemize}

\item $F$ has CBB by $K_F=\sup_B\{ |D f|^2 {\displaystyle+} Kf^2\}$, 

\end{enumerate}
\end{theorem}
\begin{remark}\hspace{0mm}
	\begin{enumerate}
		\item To illustrate the boundary condition in (1)(c) we consider the smooth context, i.e. $B$ and $F$ are Riemannian manifolds with convex boundary and sectional curvature bounded from below by $K$ and $K_\sF$ respectively. Let us assume that $f>0$ everywhere on $\partial B$. By the doubling space theorem of Perelman $B^\dagger$ is $CBB(K)$. The condition $(\dagger)$ for $f$ 
		means that the derivative in direction of the inward normal vector on $\partial B\backslash f^{-1}(\{0\})$ is negative. 
		\item If one assumes $\partial B\subset f^{-1}(\{0\})$ then the condition (1)(c) in the theorem becomes redundant.  Since $f$ is $Kf$-concave, in fact by the maximum principle it follows $\partial B= f^{-1}(\{0\})$. 
		\item The theorem has a counterpart in the context of metric spaces with curvature bounded from above. 
	\end{enumerate}
\end{remark}
\begin{example}
Take $B=[1,\infty)$ and $f(r)=r$. Then $B\times_r \mathbb{S}^1\simeq \R^2 \backslash B_1(0).$
The theorem  does not apply for any choice of $K$ since $f$ does not
satisfy $(\dagger)$; speciﬁcally, the tautological extension of $f$ to the double of
$B$ has the form of an absolute value function and cannot be $fK$-concave.
\end{example}
\begin{example}
Theorem \ref{th:albi1} is  very effective  to produce examples of Alexandrov spaces. It applies directly to all of the cones that were presented above. For instance, the Euclidean cone over some intrinsic metric space (complete, locally compact) has CBB below $0$ if and only if $F$ has CBB by $1$,  or the parabolic cone over $F$ has CBB by $-1$ if and only if $F$ has CBB by $0$. 
\end{example}
\begin{example}
A class of examples arises from the function $d_{\partial B}$, the distance in an Alexandrov space $B$ to the boundary of $\partial B$ . If $B$ has CBB by $1$ and if we set $f=\sin d_{\partial B}$, then $B\times_f F$ has CBB by $1$ provided $F$ has curvature bounded from below by $1$, or provided $F$ is an interval of length less than $\pi$. This construction generalizes the example (4) in  Example \ref{examples1} above. Explicitly, we can choose $B=\overline{\mathbb S^2_+}$ and $f=\sin d_{\partial B}$, and for $F$ we  pick any Alexandrov space with CBB by $1$. 
\end{example}
\begin{example}
The tangent cones in a warped product $B\times_f F$ are given as warped products themself. More precisely
$$T_{(s,x)}(B\times_f F)= \begin{cases} T_sB \times_{D_sf} F & \mbox{ if } f(s)=0, \\
	T_sF\times T_x F& \mbox{ if } f(s)>0.
	\end{cases}.
	$$
	We can consider $B=\overline{\mathbb S^2_+}$ and $f$ as in the previous example and  choose $\R P^2$ with the standard round metric. The warped product $B\times_f \R P^2$ is a nonsmooth Alexandrov space with CBB by $1$. Short segments in $\partial B=\mathbb S^1$ embed as minimal geodesics in this warped product via $s\in \partial B\mapsto (s,x)$. Since $f(s)$ vanishes for all $s\in \partial B$, the tangent cones along such geodesics are $\R \times [0,\infty)\times_g \R P^2\simeq \R\times C(\R P^2)$ where $g: (x,y)\in \R\times [0,\infty)\mapsto y\in [0,\infty)$.
	For instance, this shows that tangent cones in the interior of geodesics not necessarily are isometric to Euclidean spaces. 
\end{example}


%

\subsection{Curvature-dimension conditions} We review briefly the theory of metric measure spaces satisfying a synthetic Ricci curvature bound. There is a rich literature on the subject. For more details we refer the reader, for instance,  to \cite{viltot, stugeo1, stugeo2, agsriemannian, giglisplitting, ohtmea}.

For $\kappa\in \mathbb{R}$ let $\sin_{\kappa}:[0,\infty)\rightarrow \mathbb{R}$ be the solution of 
$
v''+\kappa v=0, \ v(0)=0 \ \ \& \ \ v'(0)=1.
$
For $K\in \mathbb{R}$, $N\in (0,\infty)$ and $\theta> 0$  we define the \textit{distortion coefficient} as
$$
	t\in [0,1]\mapsto \sigma_{K,N}^{(t)}(\theta)=\begin{cases}
		\frac{\sin_{K/N}(t\theta)}{\sin_{K/N}(\theta)}\ &\mbox{ if } \theta\in (0,\pi_{K/N}),\\
		\infty\ & \ \mbox{otherwise}.
	\end{cases}
$$
One sets $\sigma_{K,N}^{(t)}(0)=t$.
For $K\in \mathbb{R}$, $N\in [1,\infty)$ and $\theta\geq 0$ the \textit{modified distortion coefficient} is defined as
$$
	t\in [0,1]\mapsto \tau_{K,N}^{(t)}(\theta)=\begin{cases}
		\theta\cdot\infty \ & \mbox{ if }K>0\mbox{ and }N=1,\\
		t^{\frac{1}{N}}\left[\sigma_{K,N-1}^{(t)}(\theta)\right]^{1-\frac{1}{N}}\ & \mbox{ otherwise}.
\end{cases}$$

Let $(X,d)$ be a complete separable metric space equipped with a locally finite Borel measure $\m$. We call the triple $(X,d,\m)$ a metric measure space. 

%
%
%

We denote by $\mathcal{P}(X)$ the set of Borel probability measures on $X$ and by $\mathcal{P}_\mathrm{c}(X)$ the subset of compactly supported ones. Given a metric measure space $(X,\de,\m)$, we also define the subset $\mathcal{P}^{ac}(X)$  of $\m$-absolutely continuous probability measures and its subset $\mathcal{P}^{ac}_\mathrm{c}(X)$  of measures with compact support. For $p\geq1$, we  define the set $\mathcal{P}_p(X)$ of probability measures $\mu\in \mathcal{P}(X)$ with $\int_X\de(x,x_0)^p\,\dint\mu(x)<\infty$ for some $x_0\in X$. Similarly we define $\mathcal P^{ac}_p(X)$  and $\mathcal P^{ac}_{p,\textnormal c}(X)$.

Given $\mu,\nu\in\mathcal{P}(X)$, let $\Pi(\mu,\nu)$ be the set of all their \emph{couplings}, i.e. measures $\pi\in\mathcal{P}(X^2)$ such that $\pi(\cdot\times X)=\mu$ and $\pi(X\times\cdot)=\nu$. Alternatively, given a Borel map $T:X\rightarrow Y$ from a measure space $X$ to a measurable space $Y$, define the \emph{pushforward} $T_\sharp \mu$ of $\mu\in\mathcal{P}(X)$ by $T$ as a probability measure on $Y$ by the formula
$$T_\sharp \mu(E):=\mu(T^{-1}(E)), \ \ E\subset X \mbox{ measurable}.$$
Then $\pi\in\Pi(\mu,\nu)$ if and only if $(P_1)_\sharp \pi=\mu$ and $(P_2)_\sharp \pi=\nu$, where $P_1$ and $P_2$ are the projections onto the first and second factors. 

For $\mu,\nu\in\mathcal{P}_p(X)$ the \emph{$p$-Wasserstein distance} $W_p:\mathcal{P}_p(X)\times\mathcal{P}_p(X)\rightarrow\R$ is defined as
\begin{align}\label{id:pwasserstein}W_p(\mu,\nu):=\left(\inf_{\pi\in\Pi(\mu,\nu)}\int_{X^2}\de(x,y)^p\,\dint\pi(x,y)\right)^\frac{1}{p}.\end{align}
If $\pi\in\Pi(\mu,\nu)$ is a minimizer in  \eqref{id:pwasserstein}, we say that $\pi$ is $p$-optimal and write $\pi\in\Pi^{opt}_p(\mu,\nu)$. The function $W_p$ is a (finite) distance on $\mathcal{P}_p(X)$ and turns the latter into a complete metric space \cite{viltot}. Moreover, if $(X,\de)$ is geodesic, then also $(\mathcal{P}_p(X),W_p)$ is geodesic \cite{viltot}.

%

A dynamical optimal coupling is a probability measure $\Pi\in \mathcal P(\mathcal G(X))$  such that $t\in [0,1]\mapsto (e_t)_{\#}\Pi$ is a $W_2$-geodesic in $\mathcal P^2(X)$. 

The set of  dynamical optimal couplings $\Pi\in \mathcal P(\mathcal G^{}(X))$ between $\mu_0,\mu_1\in \mathcal P^2(X)$ is denoted with ${opt}Geo(\mu_0,\mu_1)$. 

A metric measure space $(X,d,\m)$ is called \textit{essentially nonbranching} if for any pair $\mu_0,\mu_1\in \mathcal P^2(X,\m)$ any $\Pi\in {opt}Geo(\mu_0,\mu_1)$ is concentrated on a set of nonbranching geodesics.
\smallskip

{The \textit{$N$-Renyi entropy} is defined by
	$$
	S_N(\cdot|\m):\mathcal{P}^2_b(X)\rightarrow (-\infty,0],\ \ S_N(\mu|\m)=\begin{cases}-\int_X \rho^{1-\frac{1}{N}}d\!\m& \ \mbox{ if $\mu=\rho\m$,  }\smallskip\\
		0&\ \mbox{ otherwise}.
	\end{cases}
	$$}

\begin{definition}[\cite{stugeo2, lottvillani}]\label{def:cd}
	A metric measure space $(X,d,\m)$ satisfies the \textit{curvature-dimension condition} $\CD(K,N)$ for $K\in \mathbb{R}$, $N\in [1,\infty)$ if for all $\mu_0,\mu_1\in \mathcal{P}_b(X,\m)$ 
	there exists a probability measure $\Pi\in \mathcal P(\mbox{Geo}(X))$ such that $(e_i)_\#\Pi= \mu_i$ for $i=0,1$ and $\forall t\in (0,1)$
	%
	\begin{align}\label{ineq:cd}
	S_N(\mu_t|\m)\leq\! -\!\!\int \left[\tau_{K,N}^{(1-t)}(\theta)\rho_0(e_0(\gamma))^{-\frac{1}{N}}+\tau_{K,N}^{(t)}(\theta)\rho_1(e_1(\gamma))^{-\frac{1}{N}}\right]d\Pi(\gamma)
	\end{align}
	where $\mu_i=\rho_id\m$, $i=0,1$,  $\mu_t= (e_t)_\#\Pi$, and $\theta= L(\gamma)$.
	%
\end{definition}
\begin{definition}
	A metric measure space $(X,\de, \m_\sX)$ satisfies the Riemannian curvature-dimension condition  $
	\RCD(K,N)$ if  the curvature-dimension condition $\CD(K,N)$ is coupled with the assumption that the space $X$ is infinitesimally Hilbertian, i.e. the space of Sobolev function $W^{1,2}(X)$ on $X$ is Hilbert space. 
	\end{definition}\begin{remark}
The Riemannian curvature-condition $\RCD(K, \infty)$,  that we do not  present here, was introduced in \cite{agsriemannian} coupling the curvature-dimension condition $\CD(K, \infty)$, previously proposed in \cite{stugeo1, stugeo2} and independently in \cite{lottvillani}, with the linearity of
the heat flow.
The finite dimensional refinements subsequently led to the notions of  \linebreak[4]$\RCD(K,N)$  spaces,
corresponding to $\CD(K, N)$  coupled with the infinitesimally Hilbertian condition. The
class $\RCD(K,N)$ as well as the notion of infinitesimally Hilbertian was proposed in \cite{giglistructure}. 
\end{remark}
\begin{remark}
 In \cite{Kap-Ket-18} it was shown that the Riemannian curvature-dimension condition is equivalent to the curvature-dimension condition coupled with the assumption that at $\m_\sX$-a.e. point $x\in X$ there exists an Euclidean tangent cone.  
\end{remark}
A weaker version than the $\CD$ condition is obtained by considering a convexity property just along Wasserstein geodesics  where one measure is a Dirac measure. 

\begin{definition}\cite{ohtmea, stugeo2}
    A metric measure space satisfies the measure contraction property $\MCP(K,N)$ for $K\in\R$ and $N\geq1$ if for any $o\in\supp\,\m$ and $\mu_0\in\mathcal{P}(X)$ of the form $\mu_0=\frac{1}{\m(A)}\m|_A$ for some Borel set $A\subset X$ with $0<\m(A)<+\infty$ there exists $\Pi\in \mathcal P(\Geo(X))$ such that
    \begin{align}
        \frac{1}{\m(A)}\m\geq(e_t)_\sharp\left(\tau_{K,N}^{(1-t)}(\de(\gamma_0,\gamma_1))^N\Pi(\de\gamma)\right)\ \ \ \forall t\in[0,1].
    \end{align}
\end{definition}
\subsubsection{$N$-warped products and the curvature-dimension condition}
Let $F$ be a metric measure space, and let $B$ be a $d$-dimensional  Alexandrov  space together with a continuous function $f: B\rightarrow [0,\infty)$. We want to introduce a measure on $B\times_f F$ such that the metric warped product becomes a metric measure space.  For this we will follow the construction of Section \ref{subsubsec:Nwp} for smooth spaces.  

We define the $N$-warped product between $B, f$ and $F$ as the metric measure space given by 
$$(B\times_f F, f^N \mathcal H^d_{\sB}\otimes \m_{\sF})=: B\times^\sN_f F$$
where $\mathcal H^d_\sB$ is the $d$-dimensional Hausdorff measure of $B$. We will also write frequently $f^N\de\mathcal H^d_{\sB}\otimes \de\m_\sF= :\m^\sN$.  When $B$ is smooth then this measure coincides with the $N$-warped product measure defined in Section \ref{subsubsec:Nwp}. 

For the next theorem we will consider a weighted Riemannian manifold $(F, g_\sF, \m_\sF)$, a Riemannian manifold $B$ with boundary and a smooth function $f: B\rightarrow \R$ such that $\partial B\subset f^{-1}(\{0\})$.  Since $f$ is allowed to vanisch, the resulting space will not be a smooth manifold in general. In particular, one will not be able to define lower Ricci curvature bounds in general. But for the curvature dimension condition there is the following theorem.
\begin{theorem}[\cite{ketterer1}]
Assume $(F, d_F,\m_F)$ is a weighted Riemannian manifold $(F, e^{-\Phi})$. Let $B$ be a smooth manifold with boundary  and let $f: B\rightarrow [0, \infty)$ be smooth.
\smallskip
\\
Then the {$N$-warped product} of $B, F$, $f$ satisfies $\CD(KN,N+1)$ { if }
\begin{enumerate} 
\item \begin{enumerate}\item$B$ has CBB by $K$,
	\item $f$ is $fK$-concave; 
	\end{enumerate}
\item $(F, e^{-\Phi} \vol_F)$ satisfies $\CD(K_F (N-1), N)$ with \begin{center}$K_F=\sup_I(Kf^2  + (f')^2)$.\end{center} 
\end{enumerate}
\end{theorem}
For the case that the base is a $1$-dimension interval we have the following  more general theorem where $f$ is a smooth function and $F$ is  $\RCD$.
\begin{theorem}[Ketterer]
The $N$-warped product of $I, F$ and $f$ satisfies $\RCD(K N, N+1)$ if
\begin{enumerate} 
\item $ f'' + K f \leq 0$. 
\item $(F, e^{-\Phi} \vol_F)$ satisfies $\RCD(K_F (N-1), N)$ with \begin{center}$K_F=\sup_I(Kf^2  + (f')^2)$. \end{center}
\end{enumerate}
\end{theorem}
The proof of the previous theorem follows  closely \cite{ketterer2} and will appear in an upcoming publication.

Under the assumption that $f$ satisfies $f''+K f=0$, the converse theorem holds as well. 
\begin{theorem}[\cite{ketterer2}]
Let $X$ be a metric measure space. 

The Euclidean $N$-cone 
$C(X)$ satisfies $\RCD(0, N+1)$ if and only if $X$ satisfies $\RCD(N-1,N)$.
\end{theorem}
The previous theorems suggest that the following general statement should be true. We always assume that $f$ is Lipschitz. 
\begin{conjecture}
The {$N$-warped product} of $B, F$, $f$ satisfies $\CD(K(N+d-1),N+d)$ { if }
\begin{enumerate} 
\item \begin{enumerate}\item[(a)] $B$ is $CBB(K)$, 
	\item [(b)] $f$ is $fK$-concave; 
\item[(c)] If $B^\dagger$ is the result of gluing two copies of $B$ along the closure of the set of boundary points where $f$ is nonvanishing, and $f^\dagger: B^\dagger \rightarrow [0,\infty)$ is the tautological extension of $f$, then $B^\dagger$ has CBB by $K$ and $f^\dagger$ is $fK$-concave. 
	\end{enumerate}
\item $F$ is an essentially non-braching $\CD(K_F (N-1), N)$ space with \begin{center}$K_F=\sup_B(Kf^2  + |Df|^2)$.\end{center} 
\end{enumerate}Here $Df$ is the Alexandrov differential. 
\end{conjecture}

A particular case of the previous conjecture is when $B=I$ is an interval. Especially we conjecture the following statement. 
\begin{conjecture}
Let $X$ be an essentially non-branching metric measure space.  

The Euclidean $N$-cone 
$C(X)$ satisfies $\CD(0, N+1)$ if and only if $X$ satisfies $\CD(N-1,N)$.
\end{conjecture}
A result for the $\MCP$ condition was proven by Ohta in \cite{ohtpro}.
\begin{theorem}
Let $X$ be a metric measure spaces that satisfies $\MCP(N-1,N)$. Then $C(X)$ satisfies $\MCP(0,N+1)$. 
\end{theorem}
We call an $\RCD(K,N)$ space $(X,d,\m)$ non-collapsed if $N\in \N$ and $\m=\mathcal H^N$, the $N$-dimensional Hausdorff measure. \cite{GP-noncol}

If $X$ is a non-collapsed $\RCD(K,N)$ space then every tangent cone is an Euclidean
 cone $C(Y)$ over some $\RCD(N-1, N)$ space $Y$ \cite{ketterer2, DGi}. Moreover, there is a natural stratification 
$\mathcal S^0 \subset \dots \subset \mathcal S^{N-1} = X\backslash \mathcal R$
where
$$\mathcal R=\{x\in X: \exists\mbox{ a unique GH tangent cone isometric to }\R^N\}, $$
the set of regular points in $X$, and for $k\in \{0, \dots, N-1\}$
$$\mathcal S^k= \{x\in X: \mbox{ there is no GH tangent cone isometric to } Y\times \R^{k+1}\},$$
the set of singular points of order $k\in \N$. 

For a noncollapsed Gromov-Hausdorff Ricci limit of a sequence of closed Riemannian manifolds with no boundary it was shown in \cite{cheegercoldingI} that the top-dimensional singular set $\mathcal S^{N-1}\backslash \mathcal S^{N-2}$ is empty. In a general noncollapsed $RCD$ space this set may be non-empty. For instance, one can consider a Riemannian manifold $(M,g)$ with convex boundary and Ricci curvature bounded from below. For Alexandrov spaces the top dimensional singular stratum is the Alexandrov boundary.

In \cite{bns} and \cite{Kap-Mon19} two different notions of boundary have been proposed.  In \cite{bns} the authors define 
$$\partial X= \overline{\mathcal S^{N-1}\backslash \mathcal S^{N-2}}, $$
the topological closure of $\mathcal S^{N-1}\backslash \mathcal S^{N-2}$.  

In \cite{Kap-Mon19} an alternative definition of boundary has been proposed, inspired by the one adopted for Alexandrov spaces. The authors of \cite{Kap-Mon19} define the boundary of $X$ through
$$\mathcal F X= \{x\in X: \exists\mbox{ a tangent cone }C(Z_x) \mbox{ s.t. } Z_x \mbox{ has nonempty boundary}\}.$$
Recursively the definition of boundary  reduces to the case of $RCD(0,1)$ spaces that are isometric to one dimensional manifolds thanks to a classification of such spaces given in \cite{Kit-Lak}.

 \begin{theorem}[\cite{bns}] 
Let $X$ be a noncollapsed $RCD(K,N)$ space. Then either $\partial X=\emptyset$, or $\mathcal F X\neq \emptyset$ and $\mathcal FX\subset \partial X$. 
\end{theorem}
Now it is  possible to formulate a conjecture in the following form where the assumptions on the base are further relaxed. 
\begin{conjecture}
The {$N$-warped product} of $B, F$, $f$ satisfies $\CD(K(N+d-1),N+d)$ { if }
\begin{enumerate} 
\item \begin{enumerate}\item[(a)] $B$ is a non-collapsed $\RCD((d-1)K, d)$ spaces, 
	\item [(b)] $f: B\rightarrow [0, \infty)$ is $fK$-concave; 
\item[(c)] If $B^\dagger$ is the result of gluing two copies of $B$ along the closure of the set of boundary points where $f$ is nonvanishing, and $f^\dagger: B^\dagger \rightarrow [0,\infty)$ is the tautological extension of $f$, then $B^\dagger$  is a noncollapsed $\RCD((d-1),K,d)$ spaces and $f^\dagger$ is $fK$-concave. 
	\end{enumerate}
\item $F$  is essentially non-braching and satisfies $\CD(K_F (N-1), N)$ with \begin{center}$K_F=\sup_B(Kf^2  + |\nabla f|^2)$.\end{center} 
\end{enumerate}
Here, the boundary of $B$ is understood as $\mathcal FX$ or as the closure of the top-dimensional stratum. 
\end{conjecture}
Proposition  \ref{prop:alt}  suggests the following modification of the previous conjecture. 
\begin{conjecture}
Let $(B^d, f^{N+d} \vol_\sB)$ be a weighted Riemannian manifold, possibly with boundary, for $0\leq f \in C^\infty(B)$. 

The {$N$-warped product} of $B, F$, $f$ satisfies $\CD(K(N+d-1),N+d)$ { if }
\begin{enumerate} 
\item \begin{enumerate}\item[(a)] $(B, f^{N+d}\vol_\sB)$ is a  $\RCD((N+d-1)K, N+d)$ spaces, 
	\item [(b)] $f$ satisfies $\Delta^\sB f\leq - dKf$,
\item[(c)] $\partial B \subset f^{-1}(\{0\})$;
	\end{enumerate}
\item $F$ is essentially nonbraching and satisfies $\CD(K_F (N-1), N)$ with \begin{center}$K_F=\sup_B(Kf^2  + |\nabla f|^2)$.\end{center} 
\end{enumerate}
\end{conjecture}
One notices that $Kf$-concavity of $f$ in the previous conjecture is traded for the property that $(B, f^{N+d} \vol_\sB)$ satisfies  $\RCD((N+d-1)K, N+d)$ and $f$ satisfies $\Delta f\leq - Kd f$.
%
%
%
\section{Lorenztian length spaces}
\subsection{Warped products as Lorentzian length spaces}
Let $F$ be a metric space. 
If the base space $B$ is an interval $I\subset\R$ and $f: I\rightarrow [0, \infty)$ we call  the warped product $I\times_f F$ a generalized cone of $I$, $F$ and $f$. 

In particular, if $I=(a,b)$ with $a,b\in \R\cup\{\pm\infty\}$ and $F$ is smooth Riemannian manifold then $I\times_f F$ is also a smooth Riemannian manifold. 


In the construction of a generalized cone the interval $I$ is considered with its natural Riemannian structure, i.e. the standard Riemannian metric $g_I= (\de t)^2$. One can replace $g_I$ with the semi-Riemannian metric $-g_I$ that has  index $1$. In the following we will denote $I$ equipped with $-g_I$ by $-I$.

A smooth warped product between $-I$, $f$ and a Riemannian manifold $F$ is well defined and gives a semi-Riemannian manifold of signature  $(+, -, \dots, -)$, i.e.  a Lorentzian manifold. 
We will generalize this construction now to allow non-smooth metric spaces.  This will follow the plan in \cite{agks} (see also \cite{beran}).

Let $F$ be an intrinsic metric space and let $I\subset \R$ be an open interval. We consider $I\times F$ with the product topology induced by the product metric. Let $f: I \rightarrow (0, \infty)$ be Lipschitz continuous. 

The goal is to define  a Lorentzian length structure in the sense of \cite{kusae}. To this end we first define a class of causal curves. Let $\gamma=(\alpha, \beta):[a,b]\rightarrow I\times F$ be an admissible path, i.e. $\alpha$ and $\beta$ are rectifiable. We assume that $\gamma$ is parametrized such that  both $\alpha$ and $\beta$ are parametrized proportional to  arclength in $I$ and $F$, respectively.  In particular, the metric speed of $ \alpha$ and of $\beta$ is well-defined a.e. in $[a,b]$. 
We call $\gamma$ 
\begin{align*}
\begin{cases} timelike\\
null \\
causal 
\end{cases} \mbox{ if } \  -(\alpha')^2 + (f\circ \alpha)^2 |{\beta}'|_\sF^2\ \  \begin{cases} <0 \\
=0\\
\leq 0
\end{cases} a.e. 
\end{align*} 
The path $ \gamma$ is called future/past directed if $ \alpha$ is strictly monotonically increasing/decreasing, i.e. $ \alpha'>0$ or $ \alpha'<0$ a.e.  

Points $(s,x)$ and $(t,y)$ in $I\times F$ are chronologically related, denoted by $(s,x)\ll (t,y)$, if there exists a future directed timelike curve from $(s,x)$ to $(t,y)$. Moreover, points $(s,x)$ and $(t,y)$ are causally related, denoted with $(s,x)\leq (t,y)$, if there exists a future directed causal curve from $(s,x)$ to $(t,y)$, or $(s,x)= (t,y)$. 

The cronological and causal future of a point  $(s,x)$ are defined as \begin{center}
$I^{+}((s,x))= \{(t,y): (s,x)\ll  (t,y)\},  \ \ J^+((s,x))= \{ (t,y): (s,x) \leq (t,y)\}$\end{center}
and similarly for the chronological and causal past. 

Let $\gamma=(\alpha, \beta)\rightarrow I\times F$ be a future directed causal curve. We define the {\it length } of a causal curve $\gamma=(\alpha, \beta): [a,b]\rightarrow I\times F$ as 
\begin{align*} 
\inf \sum_{i=1}^N \left({ |\alpha(t_i)- \alpha(t_{i-1})|^2 - (f\circ \alpha(t_{i-1}))^2 \de_\sF(\beta(t_{i-1}), \beta(t_i))}\right)^{\frac{1}{2}}=L(\gamma)
\end{align*}
where we take the infimum w.r.t. all partitions $t_0=a\leq t_1\leq \dots \leq t_N=b$ of $[a,b]$.  
%

Again, from the definition of $L$ one checks the following properties:
\begin{enumerate}
\item If $c\in (a,b)$, then $L(\gamma|_{[a,c]})+ L(\gamma|_{[c,b]})= L(\gamma|_{[a,b]})$ ({\it Additivity});
\item If $\gamma$ is admissible, then $c,d\in [a,b]\mapsto L(\gamma|_{[c,d]})$ is a continuous function in $c$ and $d$;
\item $L$ is lower semi-continuous on the space of admissible paths in $B\times F$ w.r.t. to pointwise convergence of admissible paths. 
\end{enumerate} Also  the length of an admissible curve is invariant under suitable  reparametrization as it is defined via partitions. 

Let $\gamma=(\alpha, \beta):[a,b]\rightarrow I\times F$ be a future directed causal curve. In particular, $\alpha$ and $\beta$ are Lipschitz curves. Then the length satisfies
$$L(\gamma)= \int_a^b \sqrt{ (\alpha')^2 - (f\circ \alpha)^2 |\beta'|^2} \de t. $$
A time separation function $\tau$ on $I\times F$ is defined from the length $L$ as follows. For two points $(p,x)$ and $(q,y)$ such that $(p,x)\leq (q, y)$ we define the time speration function
$$\tau_{\smwp}((p,x), (q,y)):= \sup L(\gamma)$$
where the sup is w.r.t. all future directed causal curves $\gamma$ that connect  $(p,x)$ and $(q,y)$.  If there is no such curve, one sets $\tau((p,x), (q,y))=0$. $\tau_{\smwp}$  satisfies the reverse $\triangle$-inequality
$$\tau_{\smwp}((p,x), (q,y)) \geq \tau_{\smwp}((p,x), (r,z))+ \tau_{\smwp}((r,z), (q,y)).$$
For all $(p,x), (q,y), (r,z)\in I\times F$. In the following we write also $\tau_{\smwp}= \tau$. 

The time seperation function has the following properties: 
\begin{itemize}
\item $\tau(y,y')= 0$ if $y'$ and $y$ are not cuassally related, 
\item $\tau(y,y')>0$ if $y\ll y'$. 
\end{itemize}

The generalized Lorentzian cone w.r.t. $-I$, the metric space $F$ and the function $f:I\rightarrow (0, \infty)$ is 
$$(I\times F, D, \ll, \leq , \tau)=: -I\times_f F.$$
If $F$ is an intrinsic metric space, then $-I\times_f F$ is a Lorentzian pre-length space in the sense of \cite{kusae}. 

Here we collect some properties of $-I\times_f F$. We always assume that $F$ is intrinsic. 
\begin{itemize}
\item $-I\times_f F$ has the push-up property: if $(s,x)\ll (t,y)$ if and only if there exists a future directed causal curve from $(s,x)$ to $(t,y)$ with positive length. Moreover $I^{+/-}((s,x))$ is open for all $(s,t)\in I\times F$.
\item If $F$ is geodesic, then $J^{\pm}((s,x))$ is closed for all $(s,x)\in I\times F$.
\end{itemize}
The next theorem is the Lorentzian analogue of Theorem \ref{th:albi0} above.
\begin{theorem}[Kunzinger, Saemann, Graf, Alexander] Let $X$ be a strictly intrinsic metric space and let $\gamma=(\alpha, \beta)$ be a future directed, causal curve that is a maximizer w.r.t. $L$ in $-I\times_f F$ and parametrized proportional to arclength. Then
\begin{enumerate}
\item[(a)] $\beta$ is a minimizer in $F$; 
\item[(b)] (Fiber independence) $\alpha$ is independent of $F$, except for the total height, i.e. the length $L^\sF(\beta)$ of $\beta$. More precisely, if $\hat F$ is another strictly intrinsic metric space and $\hat \beta$ is a minimizing geodesic in $\bar F$ with the same length and speed as $\beta$, then $(\alpha, \hat \beta)$ is a minimizer in $-I\times_f\hat F$. 
\item[(c)] If $\gamma$ is timelike, then $\beta$ has speed $\frac{c_\gamma}{2f^2\circ \alpha}$ for a constant $c_\gamma$. 
\item[(d)] If $\gamma$ is timelike, then $-\frac{1}{2}(\alpha')^2 + \frac{1}{2f^2\circ \alpha}=E$ for a constant $E$. 
\end{enumerate}
\end{theorem}
We stated above that a generalized cone $-I\times_f F$ over an intrinsic metric space $F$ is  a Lorentzian pre-length space. If $F$ is locally compact, it was shown in \cite{agks} that $-I\times_f F$ is even a Lorentzian length spaces. 
\begin{theorem}
Let $F$ be a locally compact length space. Then $-I\times_f F$ is a strongly causal Lorentzian length space. 

Moreover, if $F$ is also complete, then $-I\times_f F$ is a Lorentzian geodesic space, i.e. any two timelike related points  $(s,x)$ and $(t,y)$ can be connected by a timelike maximal curve $\gamma$ such that $L(\gamma)=\tau((s,x), (t,y))$. 
\end{theorem}
Furthermore, it was shown in \cite{agks} that $-I\times_f F$ is globally hyperbolic if $F$ is a striclty intrinsic length space that is proper. 
\begin{corollary}
Let $-I\times_f F$ be a generalized cone, where $F$ is a locally compact and complete length space. Then $I\times_f F$ is globally hyperbolic. 
\end{corollary}
\subsection{Timlike synthetic sectional curvature bounds}
The notion of curvature bounded from below (CBB) for a Lorentzian length space mimics the notion of Alexandrov curvature bounded below for metric spaces. More precisely, a Lorentzian length space is
said to have timelike curvature bounded below by $K$ if every
point in $Y$  has a neighborhood $U$ such that:
\begin{enumerate}
\item $\tau|_{U\times U}$ is finite and continuous.
\item  Whenever $x,y\in U$ with $x\ll y$, there exists a causal curve $\alpha$ in $U$
with $L(\alpha)= \tau(x,y)$.
\item If $(x,y,z)$ is a timelike geodesic triangle in $U$ , realized by maximal
causal curves $\alpha, \beta, \gamma$ whose side lengths satisfy the appropriate size
restrictions, and if $(x', y', z')$ is a comparison triangle of $(x,y,z)$ in
$\mathbb L^2(K)$ realized by timelike geodesics $\alpha', \beta', \gamma'$, then whenever $p,q$ are
points on the sides of $(x,y,z)$ and $p', q'$ are corresponding points of
$(x', y',z')$ , we have $\tau(p,q)\leq \tau'(q',q')$. $U$ is called a comparison neighborhood.
\end{enumerate}
Here $\mathbb L^2(K)$ is the simply connected Lorentzian model space of constant sectional curvature $K$.
\begin{theorem}[Alexander-Graf-S\"amann-Kunzinger]
$-I\times_f F$ has timelike CBB by $K$ { if} 
\begin{enumerate}
\item $f$ is $-Kf$-concave, 
\item $F$ has CBB by $K_F=\sup_I\{-(f')^2 + K f^2 \}$. 
\end{enumerate}
\end{theorem}
\begin{conjecture}
$-I\times_f F$ has timelike CBB by $K$ { if and only if} 
\begin{enumerate}
\item $f$ is $-Kf$-concave, 
\item $F$ has CBB by $K_F=\sup\{-(f')^2 + K f^2 \}$. 
\end{enumerate}
\end{conjecture}
A partial result of this conjecture was proven in \cite{agks} under the assumption that $-I\times_f \mathbb M_K^2$ satisfies $CBB(K)$ where $\mathbb M_K^2$ is the $2$-dimensional, simply connected space form with curvature $K$.
\subsection{Timelike curvature-dimension conditions}

There are Lorentzian versions of the curvature-dimension condition $\CD$  and the measure contraction property $\MCP$, denoted with $T\CD$ and $T\MCP$ \cite{camolorentz}. We will not recall the details of these definitions. They  mirror the correponding ones in the positive signature case, and also characterize the  timelike Ricci curvature bounds for smooth Lorentzian spaces \cite{mccannlorentz, mondinosuhr}. For more details we refer to \cite{mccannnull, braun}.

The following results will appear in an upcoming publication.
\begin{theorem}[Calisti, K. , Saemann]Assume $F$ nonbranching. 
The  $N$-warped product of $-I, F$ and $f$ satisfies $T\MCP(-KN, N+1)$ if \begin{enumerate} 
\item $ f'' - K f \leq 0$. 
\item $(F, d_F, \m_F)$  satisfies $\CD(K_F (N-1), N)$ with $$K_F=\sup_I(Kf^2  - (f')^2).$$
\end{enumerate}
\end{theorem}
A general conjecture is as follows.

\begin{conjecture} Assume $F$ nonbranching. 
The  $N$-warped product of $-I, F$ and $f$ satisfies $T\CD(-KN, N+1)$ if \begin{enumerate} 
\item $ f'' - K f \leq 0$. 
\item $(F, d_F, \m_F)$ is non-branching and satisfies $\CD(K_F (N-1), N)$ with $$K_F=\sup_I(Kf^2  - (f')^2).$$
\end{enumerate}
\end{conjecture}

%
\bibliography{new} 
\bibliographystyle{amsalpha}
\end{document}